\documentclass[twocolumn]{autart}

\usepackage[pdftex]{graphicx}

\usepackage{amsmath}
\usepackage{amssymb}
\usepackage{amsxtra,color}

\newtheorem{lemma}{Lemma}[section]
\newtheorem{theorem}[lemma]{Theorem}
\newtheorem{proposition}[lemma]{Proposition}

\newtheorem{definition}[lemma]{Definition}

\newtheorem{remark}[lemma]{Remark}
\newcommand{\R}{\ensuremath{\mathbb{R}}}
\newcommand{\B}{\ensuremath{\mathbb{B}}}
\newcommand{\N}{\ensuremath{\mathbb{N}}}
\newcommand{\id}{\ensuremath{\mbox{id}}}
\newcommand{\dom}{\ensuremath{\mbox{dom}}}

\newcommand{\doms}{\ensuremath{\mbox{\scriptsize dom}}}

\newcommand{\A}{\ensuremath{\mathcal{A}}}
\newcommand{\K}{\ensuremath{\mathcal{K}}}
\newcommand{\KL}{\ensuremath{\mathcal{KL}}}
\newcommand{\KLL}{\ensuremath{\mathcal{KLL}}}

\newcommand{\Lip}{\ensuremath{Lip_{loc}}}

\begin{document}

\begin{frontmatter}

\title{Input-to-state stability of interconnected hybrid systems\thanksref{footnoteinfo}}

\thanks[footnoteinfo]{This paper was not presented at any IFAC
meeting. Corresponding author M.~Kosmykov. Tel. +49421-218-63746,
}

\author[SD]{Sergey~Dashkovskiy, Michael~Kosmykov}\ead{\{dsn,kosmykov\}@math.uni-bremen.de}

\address[SD]{Department of Mathematics and Computer Science, University of Bremen, 28334 Bremen, Germany}


\begin{abstract}
We consider the interconnections of arbitrary topology of
a finite number of ISS hybrid  systems and study
whether the ISS property is maintained for the overall system. We show that if the
small gain condition is satisfied, then the whole network is ISS
and show how a non-smooth ISS-Lyapunov function can be explicitly constructed in this
case.
\end{abstract}
\end{frontmatter}
\vspace{-3mm}
\section{Introduction}
\vspace{-3mm} Hybrid systems allow for a combination of continuous
and discontinuous types of behavior in one model and hence can be
used in many applications, for example in
robotics \cite{ASL93},
reset systems \cite{NZT08} or networked control systems \cite{TaN08,NeT08}. Such systems often have a large scale
interconnected structure and can be naturally modeled as
interconnected hybrid systems. In this paper our main interest is
in stability and robustness for such interconnections as these
properties are certainly of great importance for applications.
We will use the framework of input-to-state stability (ISS) that
was first introduced for continuous systems in \cite{Son89} and
then extended to other types of systems including hybrid ones, see
for example
\cite{CaT05},
 \cite{HLT08}, \cite{IPJ09} and \cite{JiW01}.
The ISS property of the interconnected systems is usually studied
using small gain conditions that take the interconnection
structure of the whole system into account. First small gain
conditions for ISS systems were introduced for the interconnections of
two continuous systems in \cite{JTP94,JMW96}. These results were
extended to arbitrary number of interconnected systems in
\cite{DRW07,DRW10,KaJ09}.\\
Interconnections of two hybrid ISS systems were considered in
\cite{LiN06}, \cite{NeT08}, e.g.. A stability condition of the
small gain type was used in \cite{NeT08} for a construction of an
ISS-Lyapunov function for their feedback connection.
Interconnection of arbitrary number of sampled-data systems that
are a special class of hybrid systems was considered in
\cite{KaJ09}. The small gain condition was given there in terms of
vector Lyapunov functions.\\
In this paper we obtain similar results
for the interconnection of arbitrary number of hybrid systems. To this end we use the methodology
recently developed in \cite{DRW07}, 
\cite{DRW10} for the investigation of stability of general networks of
ISS systems. In particular we use the small gain condition developed
in these papers and we use non-smooth ISS-Lyapunov functions in
our considerations. The main result of this paper extends the
result of \cite{NeT08} for the case of interconnection of $n\ge2$
hybrid systems and \cite{KaJ09} for general type of hybrid systems
by applying the small gain condition in the matrix form.
 Moreover,
we prove the small gain results in terms of trajectories.\\
There are different ways to introduce ISS-Lyapunov functions for
hybrid systems, see for example \cite{CaT09}, \cite{NeT08}. We show their equivalence in this paper.
Using the methods developed in
\cite{DRW10} we provide an explicit construction of an ISS-Lyapunov function
for interconnected hybrid systems. \\
The next section introduces all necessary notions and notation.
Section~\ref{sec:SGC} contains the main results and
Section~\ref{sec:concl} concludes the paper.
\vspace{-4mm}
\section{Preliminaries}\label{sec:Preliminaries}
\vspace{-4mm} Let ${\mathbb{R}_{+}}$ be the set of nonnegative
real numbers, ${\mathbb{R}}_{+}^n$ be the positive orthant $\{x
\in {\mathbb{R}}^n: x\geq 0\}$ and $\N_{+}:=\{0,1,2,\ldots\}$.
$x^T$ stands for the transposition of a vector $x \in
\mathbb{R}^n$. $\B$ is the open unit ball centered at the origin
in $\R^n$ and $\bar{\B}$ is its closure. Set $B\subset\R^n$ is called
{\it relatively  closed in} $\chi \subset\R^n$, if  $B=\bar{B}\cap
\chi$. By $\langle\cdot,\cdot\rangle$ we denote the standard
scalar product in $\R^n$. For $x,y \in {\mathbb{R}}^n$, we  write
$ x\geq y\, \Leftrightarrow\, x_i \geq y_i$;
$ x>y\Leftrightarrow\, x_i>y_i, i=1,\ldots,n$;
$x\not\geq y\Leftrightarrow\exists i{\in}\{1,\ldots,n\}:x_i<y_i$.
${M}^{n}$ denotes the $n$-fold
composition $M{\circ}{\ldots}{\circ} M$ of a map $M:\R_+^n{\to}\R_+^n$.\\
A function $\alpha: \R_+\to\R_+$ with $\alpha(0)=0$ and
$\alpha(t)>0$ for $t>0$ is called positive definite. A function
$\gamma:{\mathbb{R}}_{+}\to{\mathbb{R}}_{+}$ is said to be of
class $\K$  if it is continuous, strictly increasing and
$\gamma(0)=0$. It is of class $\K_{\infty}$ if, in addition, it is
unbounded. Note that for $\gamma \in \K_{\infty}$ the inverse
function $\gamma^{-1}\in
\K_{\infty}$ always exists.
A function $\beta:
{\mathbb{R}}_{+}{\times}{\mathbb{R}}_{+}{\to}{\mathbb{R}}_{+}$ is
said to be of class $\KL$ if, for each fixed $t$, the
$\beta(\cdot,t)\in\K$ and, for each
fixed $s$, the function $\beta(s,\cdot)$ is non-increasing and
tends to zero for $t{\to}\infty$.
A function $\beta:
{\mathbb{R}}_{+}{\times}{\mathbb{R}}_{+}{\times}{\mathbb{R}}_{+}{\to}{\mathbb{R}}_{+}$ is
said to be of class $\KLL$ if, for each fixed $r\geq 0$, $\beta(\cdot,\cdot,r)\in\KL$ and $\beta(\cdot,r,\cdot)\in\KL$.
\vspace{-0.3cm}
\subsection{Interconnected hybrid systems}
\vspace{-3mm} Consider an interconnection of $n$ hybrid subsystems
with states $x_i\in\chi_i{\subset}\R^{N_i}$, $i{=}1{,}\dots{,}n$, and
external input $u\in U{\subset}\R^{M}$. Dynamics of the $i$th
subsystem is given by
\vspace{-0.3cm}
\begin{equation}\label{is}
\begin{array}{llll}
\dot{x}_i&{=}&f_i(x_1,{\ldots},x_n,u),& \hspace{-0.1cm}(x_1,{\ldots},x_n, u)\in C_i\\
x^{+}_i&{=}&g_i(x_1,{\ldots},x_n,u),& \hspace{-0.1cm}(x_1,{\ldots},x_n, u)\in D_i
\end{array}
\vspace{-0.3cm}
\end{equation}
where  $f_i:C_i\to \R^{N_i}, g_i:D_i\to\chi_i$
and $C_i, D_i$ are the subsets of $\chi_1\times\ldots\times\chi_n\times U$.
Each subsystem is described by $(f_i,g_i,C_i,D_i,\chi_i,U)$, however
in view of stability properties we will need to restrict such
interconnections to $D_i=D_j,\,\forall i,j$, see
Remark~\ref{CD-problems} below. \\ If $(x_1,\ldots,x_n, u){\in}C_i$,
then system \eqref{is} flows continuously and the dynamics is given
by function $f_i$. If $(x_1,\ldots,x_n, u){\in}D_i$, then the system
jumps instantaneously according to function $g_i$. In points of
$C_i{\cap}D_i$
 the system may either
flow or jump, the latter only if the flowing keeps
$(x_1,\ldots,x_n, u){\in}C_i$.
Define $\chi{:=}\chi_1{\times}{\ldots}{\times}\chi_n$.
The solutions are defined on hybrid time domains \cite{GHT04}.
A subset $\R_{+}\times\N_{+}$  is called hybrid time domain denoted by {\it $\dom$} if it is given as a
union of finitely or infinitely many intervals
$[t_k,t_{k+1}]\times \{k\}$, where the numbers $0=t_0,t_1,\ldots$
form a finite or infinite, nondecreasing sequence of real numbers.
The ''last'' interval is allowed  to be of the form $[t_K,T)\times
\{k\}$ with $T$ finite or $T=+\infty$.\\
A hybrid signal is a function defined on the hybrid time domain. For
the $i$th subsystem the hybrid input
\vspace{-0.3cm}
\begin{equation}\label{eq:compisit-input}
v_i:=(x_1^T,\ldots,x_{i-1}^T,x_{i+1}^T,\ldots,x_n^T,u^T)^T,
\vspace{-0.3cm}
\end{equation}
consists of hybrid signals $u:\dom\; u{\to}
U\subset\R^{M}$, $x_j:\dom \,x_j\to \chi_j,j\neq i$ such that
$u(\cdot,k),x_j(\cdot,k)$ are Lebesgue measurable and locally
essentially bounded for each $k$.
For a signal $u{:}\dom\, u{\to} U{\subset}\R^{M}$ we define its
restriction to the interval $[(t_1,j_1),(t_2,j_2)]{\in} \dom\, u$
by\\
$u_{[(t_1,j_1),(t_2,j_2)]}(t,k){=}\hspace{-0.1cm}\left\{\hspace{-0.1cm}\begin{array}{ll}
u(t,k){,}& \mbox{if }(t_1,j_1){\leq}(t,k){\leq}(t_2,j_2),\\
0,& \mbox{otherwise,}
\end{array}\right.$
where for the elements of the hybrid time domain we define that
$(s,l)\leq(t,k)$ means $s+l\le t+k$. For convenience, we denote
$u_{(t,k)}:=u_{[(0,0),(t,k)]}$.\\
A hybrid arc of subsystem $i$ is such a hybrid signal $x_i:\dom
\,x_i\to\chi_i$, that $x_i(\cdot,k)$ is locally absolutely
continuous for each $k$. Define
$x:=(x^T_1,\ldots,x^T_n)^T\in\chi\subset\R^{N}$, $N:=\sum N_i$. A hybrid arc and a hybrid input is a solution pair
$(x_i,v_i)$ of the $i$th hybrid subsystem \eqref{is}
if\\
(i) $\dom\,x_i=\dom \,u=\dom \,x_j,j\neq i$ and
\vspace{-0.25cm}
\[(x(0,0),u(0,0))\in C_i\cup
D_i,\vspace{-0.25cm}\]
(ii) for all $k\in N_{+}$ and almost all $(t,k)\in \dom\,\: x_i$
\vspace{-0.25cm}
\begin{equation}\label{iarc}
\dot{x}_i(t,k){=}f_i(x(t,k),u(t,k)), \mbox{ if }(x(t,k),u(t,k)){\in}C_i
\vspace{-0.15cm}
\end{equation}
(iii) for all $(t,k)\in \dom\,  x_i$ such that $(t,k+1)\in \dom
\,x_i$
\vspace{-0.15cm}
\begin{equation}\label{iinput}
{x}_i(t,k{+}1){=}g_i(x(t,k),u(t,k)), \mbox{ if }(x(t,k),u(t,k)){\in}D_i.
\end{equation}
For the existence of solutions assume that the following basic
regularity conditions \cite{CaT09}, \cite{GoT06} hold :
\begin{enumerate}
\vspace{-0.2cm}
\item $\chi_i$ is open, $U$ is closed, and $C_i,D_i\subset\chi\times U$ are relatively closed in $\chi\times U$;
\item $f_i, g_i$ are continuous.
\vspace{-0.2cm}
\end{enumerate}
The supremum norm of a hybrid  signal $u$ defined on
$[(0,0),(t,k)]\in\dom \,u$ is defined by \vspace{-0.35cm}
\[\|u\|_{(t,k)}{:=}\max\hspace{-0.1cm}\left\{\hspace{-0.1cm}\mathrel{\mathop{\text{ess sup}}\limits_{\substack{(s,l)\in\doms \, u \backslash\Phi(u), \\
(s,l)\leq(t,k)}}}\hspace{-0.2cm}|u(s,l)|,\mathrel{\mathop{\text{sup}}\limits_{\substack{(s,l)\in\Phi(u),
\\ (s,l)\leq(t,k)}}}\hspace{-0.2cm}|u(s,l)|\hspace{-0.1cm}\right\}
\vspace{-0.25cm}\] and
$\Phi(u){:=}\{(s,l){\in}\dom{:}(s,l{+}1){\in}\dom \,u\}$. If
$t{+}k{\to}\infty$, then $\|u\|_{(t,k)}$ is denoted by
$\|u\|_{\infty}$.  The set of hybrid inputs in $\R^M$ with finite
$\|\cdot\|_{\infty}$ is denoted by $\mathcal{L}_{\infty}^M$. A
solution pair of hybrid system is {\it maximal} if it cannot be
extended. It is {\it complete} if its hybrid time domain is
unbounded. Let $S_u(x_0)$ be the set of all maximal solution pairs
$(x,u)$ to (\ref{ws}) with $x(0,0)=x_0$.\\
To consider interconnection \eqref{is} as one hybrid system
\vspace{-0.35cm}
\begin{equation}\label{ws}
\begin{array}{llll}
\dot{x}&=&f(x,u),& (x, u) \in C,\\
x^{+}&=&g(x,u),& (x, u) \in D,
\end{array}
\vspace{-0.35cm}
\end{equation}
with state $x$ and input $u$ defined above, it seems to be natural
to define $C:=\cap C_i$, $D:=\cup D_i$, since a jump of any
subsystem means a jump for the overall state $x$, and to define
function $f:C\to\R^{N}$ by $f{:=}(f_1^T,\ldots,f_n^T)^T$ and
function $g:D\to\chi$ as follows
$g{:=}(\widetilde{g}_1^T,\ldots,\widetilde{g}_n^T)^T$, where
\vspace{-0.4cm}
\begin{align}\label{new-jumps}
\widetilde{g}_i(x,u):=\left\{
\begin{array}{ll}
g_i(x,u), & \mbox{ if } (x,u)\in D_i,\\
x_i,& \mbox{ otherwise }.
\end{array}
\right.
\end{align}
\vspace{-0.4cm}\\
Note that the solutions of \eqref{ws} may have different hybrid
time domains than the solutions of the individual systems
\eqref{is}, see \cite{San10}. The above choice of $C$ and $D$ was
used also in \cite{San10} considering interconnections of two
hybrid systems. However this choice has certain drawbacks:  \eqref{new-jumps} rules out solutions starting in $C\cap D$
such that one subsystem jumps while another one flows, another problem is discussed in Remark~\ref{CD-problems}, see also \cite[Remark~4.3]{San10}. \\
\vspace{-0.6cm}
\begin{remark}\label{CD-problems}
Let one of the subsystems, say the $j$th one, has the property that
once being in $D_j$ it makes only jumps and never leaves $D_j$. If
$D_i{\neq} D_j$ for some $i$, then for any initial state $x(0,0)$
with $(x(0,0),u(0,0)){\in} D_j$ and $(x(0,0),u(0,0)){\in}
C_i\setminus D_i$ there exists a solution pair given by
$(x(0,k),u(0,k)),\,k{\in}\N$ with $x_i(0,k){=}x_i(0,0),\,\forall k$,
i.e., a solution with the "frozen" $x_i$. This follows from
\eqref{new-jumps}: being in $C_i$ we have $\tilde
g_i{=}\id$. This in particular shows that even in case of a zero
input signal there is a solution that will never become "small",
contradicting the ISS or the AG property (defined below). For this reason we require
in Section~\ref{sec:SGC} that the jump sets $D_i$ coincide for all
subsystems. This requirement implies that the subsystems can jump
simultaneously only. This restricts the class of interconnected
systems considered in this paper.
\end{remark}
\subsection{Input-to-state stability and Lyapunov functions}
\vspace{-3mm} To study stability of the interconnected hybrid
systems we use the notion of input-to-state stability (ISS)
\cite{CaT09}: \vspace{-2mm}
\begin{definition}
The $i$th subsystem (\ref{is}) is called ISS, if there exist
$\beta_i\in\KLL$, $\gamma_{ij}$, $\gamma_i \in
\K_\infty\cup\{$0$\}$ such that for all initial values $x_{i0}$
each solution pair $(x_i,v_i)\in S_{v_i}(x_{i0})$ with $v_i$ from
\eqref{eq:compisit-input} satisfies $\forall (t,k)\in\dom\,x_i$
the following:
\vspace{-2mm}
\begin{equation}\label{iISS}
|\hspace{-0.025cm}x_i(\hspace{-0.025cm}t\hspace{-0.025cm},\hspace{-0.025cm}k\hspace{-0.025cm})\hspace{-0.025cm}|\hspace{-0.025cm}{\leq}\hspace{-0.025cm}\max\{\hspace{-0.025cm}\beta_i(\hspace{-0.025cm}|\hspace{-0.025cm}x_{i0}\hspace{-0.025cm}|\hspace{-0.025cm},\hspace{-0.025cm}t\hspace{-0.025cm},\hspace{-0.025cm}k\hspace{-0.025cm}){,}\hspace{-0.025cm}\max\limits_{j{,}j{\neq}
i}\hspace{-0.025cm}\gamma_{ij}({\hspace{-0.025cm}\|\hspace{-0.025cm}x_{j}\hspace{-0.025cm}\|_{\hspace{-0.025cm}(t,k)\hspace{-0.025cm}}}){,}\hspace{-0.025cm}\gamma_i\hspace{-0.025cm}({\hspace{-0.025cm}\|u\|}_{\hspace{-0.025cm}(t,k)\hspace{-0.025cm}}\hspace{-0.025cm})\hspace{-0.025cm}\}.
\vspace{-1mm}
\end{equation}
Functions $\gamma_{ij},\gamma_i$ are called ISS
nonlinear gains.
\end{definition}
\vspace{-0.2cm}
We borrow also the following stability notions from \cite{CaT09} that
will be used in the next section to prove one of the main results (Theorem~\ref{the:main_traj}):
%
\vspace{-0.2cm}
\begin{definition}
System (\ref{is}) is called {\it $0$-input pre-stable}, if for
each $\epsilon_i>0$ there exists $\delta_i>0$ such that each solution
pair $(x_i,0)\in S_{v_i}(x_{i0})$ with $|x_{i0}|\leq\delta$ satisfies
$|x_i(t,k)|\leq\epsilon_i$ for all $(t,k)\in\dom\,x_i$.
\end{definition}
\vspace{-0.3cm}
\begin{definition}
System (\ref{is}) is called {\it globally pre-stable} (pre-GS), if
 $\exists\sigma_i,\hat{\gamma}_{ij},\hat{\gamma}_i{\in}
\K_\infty{\cup}\{$0$\}$ such that for all initial values $x_{i0}$
each solution pair $(x_i,v_i)\in S_{v_i}(x_{i0})$ satisfies $\forall (t,k)\in\dom\,x_i$
the following:
\vspace{-0.2cm}
\begin{equation}\label{iGS_max}
|x_i(t,k)|\hspace{-0.025cm}{\leq}\hspace{-0.025cm}\max\{\sigma_i(\hspace{-0.025cm}|x_{i0}|\hspace{-0.025cm})\hspace{-0.025cm},\hspace{-0.025cm}\max\limits_{j,j\neq
i}\hat{\gamma}_{ij}({\|x_{j}\|_{(t,k)}})\hspace{-0.025cm},\hspace{-0.025cm}\hat{\gamma}_i({\|u\|}_{(t,k)})\hspace{-0.025cm}\}
\end{equation}
\end{definition}
\vspace{-0.3cm}
\begin{remark}\label{remark:equiv}
Note that pre-GS follows from ISS by taking $\sigma_i(|x_{i0}|):=\beta_i(|x_{i0}|,0,0)$ and $0$-input pre-stability follows from pre-GS by considering $x_j=0,u=0$.
\end{remark}\vspace{-0.3cm}
\begin{definition}
System (\ref{is}) has the {\it asymptotic gain property} (AG), if
there exist $\widetilde{\gamma}_{ij},\widetilde{\gamma}_i{\in}
\K_\infty{\cup}\{$0$\}$ such that for all initial values $x_{i0}$
all solution pairs $(x_i,v_i){\in}S_{v_i}(x_{i0})$ are bounded
and, if complete, then satisfy
\vspace{-0.3cm}
\begin{equation}\label{iAS_max}
\limsup\limits_{\substack{(t,k)\in\doms\, x,\\
t+k\to\infty}}\hspace{-0.25cm}|x_i(t\hspace{-0.04cm},\hspace{-0.04cm}k)|\hspace{-0.02cm}{\leq}\hspace{-0.05cm}\max\hspace{-0.04cm}\{\max\limits_{j,j{\neq}
i}\hspace{-0.04cm}\widetilde{\gamma}_{ij}(\hspace{-0.04cm}{\|x_{j}\|_{\infty}}\hspace{-0.04cm})\hspace{-0.04cm}{,}\hspace{-0.04cm}\widetilde{\gamma}_i({\hspace{-0.04cm}\|u\hspace{-0.02cm}\|}_{\infty})\hspace{-0.04cm}\}.
\vspace{-0.35cm}
\end{equation}
\vspace{-0.35cm}
\end{definition}
\vspace{-0.25cm}
In Theorem~3.1 in \cite{CaT09} the following relation between ISS and AG with $0$-input pre-stability was proved.
\vspace{-0.3cm}
\begin{theorem}\label{thm:ISS-equiv}
Let the set $\{f_i(x,u):u\in U\cap\epsilon\bar{\B}\}$ be convex
$\forall\varepsilon>0$ and for any $x\in\chi$. Then \eqref{is} is
ISS if and only if it has the AG property and it is $0$-input
pre-stable.
\end{theorem}
\vspace{-0.3cm}
A common alternative to prove ISS is to use ISS-Lyapunov functions as defined below.
We consider locally Lipschitz continuous functions $V_i:
\chi_i\to\R_{+}$ that are differentiable almost everywhere by the
Rademacher's theorem. The set of such functions we denote by
$\Lip$. In points where such a function is not differentiable we
use the notion of Clarke's generalized gradient, see \cite{CLS98},
\cite{DRW10}. The set \vspace{-0.3cm}
\begin{equation}\label{def:clarke_gradient}
{\partial}V_i\hspace{-0.025cm}(\hspace{-0.025cm}x_i\hspace{-0.025cm})\hspace{-0.025cm}{=}\mbox{conv}\{\hspace{-0.025cm}\zeta_i{\in} \R^{n_i}\hspace{-0.025cm}{:}\hspace{-0.025cm}{\exists} x_{i}^p\hspace{-0.05cm}{\to} x_i{,}{\exists}{\nabla} V_i(\hspace{-0.025cm}x_{i}^p\hspace{-0.025cm}) \,\mbox{and}\,\hspace{-0.025cm}\nabla V_i(\hspace{-0.025cm}x_{i}^p\hspace{-0.025cm})\hspace{-0.025cm}{\to}\zeta_i\hspace{-0.025cm}\}
\end{equation}
is called {\it Clarke's generalized gradient of $V_i$ at
$x_i\in\chi_i$}. If $V_i$ is differentiable at some point, then
${\partial}V_i\hspace{-0.025cm}(\hspace{-0.025cm}x_i\hspace{-0.025cm})$
coincides with the usual gradient at this point. \vspace{-0.3cm}
\begin{definition}\label{def:lyap_v_i}
Function $V_i{:}\chi_i{\to}\R_+$, $V_i{\in}\Lip$
 is
called {\it an ISS-Lyapunov function for system (\ref{is})} if\\
1) There exist functions $\psi_{i1}, \psi_{i2} \in \K_{\infty}$ s.t.:
\vspace{-0.3cm}
\begin{equation}\label{is cond 1 sum}
\psi_{i1}(|x_i|)\leq V_i(x_i)\leq  \psi_{i2}(|x_i|) \mbox{ for any } x_i\in \chi_i.
\vspace{-0.3cm}
\end{equation}
2) There exist continuous, proper, positive definite functions $V_j{:}\chi_j{\to} \R$, $V_j{\in}\Lip$, $j{\in}\{1,\ldots,n\}{\setminus}\{i\}$, functions
$\gamma_{ij}, \gamma_i {\in} \K_{\infty}$ and
continuous, positive definite functions $\alpha_i$, $\lambda_i$,
with $\lambda_i(s){<}s$ for all $s{>}0$ such that for all $(x,u){\in}
C_i$
\vspace{-0.35cm}
\begin{equation}\label{is cond 2.1}
\begin{array}{l}
V_i(x_i)\geq \max\{\max\limits_{j,j\neq i}\{\gamma_{ij}(V_j(x_j))\},\gamma_i(|u|)\}\Rightarrow\\
\forall \zeta_i\in\partial V_i(x_i): \left\langle \zeta_i,f_i(x,u)\right\rangle\leq -\alpha_i(V_i(x_i))
\end{array}
\end{equation}
and for all $(x,u)\in D_i$
\vspace{-0.35cm}
\begin{equation}\label{is cond 2.2}
\begin{array}{l}
V_i(g_i(x,u))\hspace{-0.05cm}{\leq}\hspace{-0.05cm}\max\{\hspace{-0.05cm}\lambda_i(V_i(x_i)){,}\max\limits_{j,j\neq i}\{\gamma_{ij}(V_j(x_j))\}{,}\gamma_i(|u|)\hspace{-0.05cm}\}{.}
\end{array}
\vspace{-0.1cm}
\end{equation}
Functions $\gamma_{ij}$, $\gamma_{i}$ are called ISS Lyapunov
gains corresponding to
 the inputs $x_j$ and
$u$ respectively.
\end{definition}
\vspace{-0.3cm} Note that this definition is different from the
definition of an ISS Lyapunov function used in \cite{CaT09}.
 The equivalence between their
existence for \eqref{is} is shown in Appendix, Section~\ref{section:equivalent_formulations}.
Note also that
$\gamma_{ij}$ are taken the same in (\ref{is cond 2.1}) and
(\ref{is cond 2.2}).
This can be always achieved by taking
the maximums of separately obtained $\gamma_{ij}$'s for the continuous and
discrete dynamics.
 If $V_i$ is differentiable at $x_i$, then \eqref{is cond 2.1} can be
written as
\vspace{-0.3cm}
\[V_i(x_i){\geq} \max\{\max\limits_j\{\gamma_{ij}(V_j(x_j))\},\gamma_i(|u|)\}\Rightarrow\]
\vspace{-0.5cm}
\[{\nabla} V_i(x_i)\cdot f_i(x,u){\leq} {-}\alpha_i(V_i(x_i)), (x,u)\in C_i.\vspace{-0.1cm}\]
Relations between the  existence of a smooth ISS-Lyapunov function
and the ISS property for hybrid systems were discussed in
\cite{CaT09}. Proposition~2.7 in \cite{CaT09} shows
that if a hybrid system has an ISS-Lyapunov function, then it is
ISS. Example~{3.4} in \cite{CaT09} shows that the converse is in
general not true. In \cite[Theorem~3.1]{CaT09} it was proved that
if \eqref{ws} is ISS with $f$ such that the set $\{f(x,u){:}u{\in}
U{\cap}\epsilon\bar{\B}\}$ is convex ${\forall}\varepsilon{>}0$ and for
any $x\in\chi$, then it has an ISS-Lyapunov function.
Usually Lyapunov function is required to be smooth, but smoothness can be relaxed to
locally Lipschitzness as shown below.
\vspace{-0.25cm}
\begin{proposition}\label{prop:Lyapunov for ISS}
If system (\ref{is}) has a locally Lipschitz continuous
ISS-Lyapunov function, then it is ISS.
\end{proposition}
\vspace{-0.35cm} \textit{\textbf{Sketch of proof.}} The proof of
 \cite[Proposition~2.7]{CaT09} stated with
$\alpha_i{\in}\K_{\infty}$ works without change if $\alpha_i$ is
continuous and positive definite. As well this proof can be
extended to the nonsmooth $V_i$ using the Clarke's generalized
derivative. The assertion of the proposition follows then from
this extension and Proposition~\ref{prop:lyap_equiv} from
Section~\ref{section:equivalent_formulations} in Appendix. \hfill $\Box$\\
Note that ISS of all subsystems does no guarantee ISS of their interconnection \cite{DRW07}.
In the following
section we introduce conditions that guarantee stability for
interconnections of ISS hybrid systems.
\section{Main results}\label{sec:SGC}
\vspace{-4mm} The main question of this paper is whether the
interconnection
\eqref{ws} of the ISS subsystems \eqref{is} is ISS. 
To study this question we collect the gains $\gamma_{ij}$ of the subsystems in the matrix $\Gamma=(\gamma_{ij})_{n\times n}$,
$i,j{=}1,\dots,n$ denoting $\gamma_{ii}\equiv~0$,\, $i=1,\dots,n$,
for completeness \cite{DRW07,Rue07}. The matrix $\Gamma$ describes the interconnection topology of the whole network and
contains the information about the mutual influence between the
subsystems. We also introduce the following gain operator
$\Gamma_{\max}:\R_{+}^n\to\R_{+}^n$, see \cite{DRW07,Rue07,KaJ09}:
\vspace{-0.2cm}
\begin{equation}\label{operator_gamma}
\Gamma_{\max}(s){:=}\left(\hspace{-0.1cm}\begin{array}{c}
\max\{\gamma_{12}(s_2),\ldots,\gamma_{1n}(s_n)\}\\
\vdots\\
\max\{\gamma_{n1}(s_1),\ldots,\gamma_{n,n-1}(s_{n-1})\}
\end{array}\hspace{-0.1cm}\right){.}
\end{equation}
We define the small gain condition as follows:
\begin{equation}\label{eq:SGC_max}
\Gamma_{\max}(s) \not\geq s,\quad \forall s\in\R^n_+, \; s\neq 0.
\end{equation}
This condition was introduced and
studied in
\cite{DRW07} and \cite{Rue07}. Furthermore, in \cite{DRW07} it was shown that \eqref{eq:SGC_max} is equivalent to the so-called cycle condition \cite{KaJ09}. 
 We will see that condition \eqref{eq:SGC_max} guarantees stability of the network.
\vspace{-0.3cm}
\subsection{Small gain theorems in terms of trajectories}
\vspace{-0.2cm}
The following small gain theorems extend Theorem~4.1 and Theorem~4.2 in \cite{DRW07} to the case of hybrid systems.
\vspace{-0.3cm}
\begin{theorem}\label{sgc_gs}
Consider system (\ref{ws}) and assume that all its subsystems are pre-GS. If $\Gamma_{\max}$ defined in
\eqref{operator_gamma} with $\gamma_{ij} = \hat{\gamma}_{ij}$ satisfies  \eqref{eq:SGC_max}, then (\ref{ws}) is pre-GS, i.e. for some  $\sigma,\hat{\gamma}\in\K_\infty\cup\{$0$\}$ and for all $(t,k)\in\dom\,x$
\vspace{-0.2cm}
\begin{equation}\label{ws_GS_max}
|x(t,k)|\leq\max\{\sigma(|x_0|),\hat{\gamma}({\|u\|}_{(t,k)})\}.
\end{equation}
\vspace{-0.25cm}
\end{theorem}
\vspace{-0.25cm}
\vspace{-0.4cm}
\begin{theorem}\label{sgc_ag}
Consider the interconnected system (\ref{ws}) with $D_i=D,
i{=}1{,}{\dots}{,}n$. Assume that each subsystem \eqref{is} has the AG
property and that solutions of the system \eqref{ws} exist, are
bounded and some of them are complete. If $\Gamma_{\max}$ defined
by \eqref{operator_gamma} with $\gamma_{ij} = \widetilde{\gamma}_{ij}$ satisfies \eqref{eq:SGC_max} then system
(\ref{ws}) satisfies the AG property. In particular any complete
solution with some $\widetilde{\gamma}\in\K_\infty\cup\{$0$\}$ satisfies
\vspace{-0.2cm}
\begin{equation}\label{ws_AG}
\limsup\limits_{(t,k)\in\doms \, x,t+k\to\infty}|x(t,k)|\leq\widetilde{\gamma}({\|u\|}_{\infty}).
\vspace{-0.4cm}
\end{equation}
\end{theorem}
\vspace{-0.35cm}
Note that if all solutions of \eqref{ws} are not
complete, then \eqref{ws} is AG by definition.  See
Appendix~\ref{appendix:sgc_gs}, \ref{appendix:sgc_ag} for the
proofs of Theorem~\ref{sgc_gs} and \ref{sgc_ag}. \vspace{-0.25cm}
\begin{remark}
The existence and boundedness of solutions of \eqref{ws} is essential, otherwise
the assertion is not true, see Example~14 in \cite{DRW07}.
\end{remark}
\vspace{-0.25cm} The following theorem extends \cite[Theorem~1]{LiN06} showing ISS of
an interconnection of two ISS hybrid systems under the small gain condition. Here we show that the same holds for
arbitrary finite number of hybrid systems.
\vspace{-0.25cm}
\begin{theorem}\label{the:main_traj}
Consider the interconnected system (\ref{ws}) with $D_i=D,
i=1,\dots,n$. Assume that the set $\{f(x,u):u\in
U\cap\epsilon\bar{\B}\}$ is convex for each $x\in\chi,
\epsilon>0$. If all subsystems in (\ref{is}) are ISS and
$\Gamma_{\max}$ defined in \eqref{operator_gamma} satisfies
\eqref{eq:SGC_max}, then (\ref{ws}) is ISS, i.e. for $(t,k){\in}\dom\,x$
\vspace{-0.2cm}
\begin{equation}\label{wISS}
|x(t,k)|{\leq}\max\{\beta(|x_0|,t,k){,}\gamma({\|u\|}_{(t,k)})\}
\vspace{-0.2cm}
\end{equation}
\vspace{-0.2cm}
holds for some $\beta\in\KLL$ and $\gamma\in\K_\infty\cup\{$0$\}$.\\
\end{theorem}
\vspace{-0.4cm} \textit{\textbf{Sketch of proof.}} The idea
follows from the proof of a similar theorem for continuous systems
in \cite{DRW07}. We describe it briefly: By
Remark~\ref{remark:equiv} and Theorem~\ref{thm:ISS-equiv}, since
each subsystem is ISS, they are pre-GS and have the AG property.
By Theorem~\ref{sgc_gs} and Theorem~\ref{sgc_ag} the whole
interconnection \eqref{ws}  is pre-GS and has the AG property.
From global pre-stability of \eqref{ws}, $0$-input pre-stability
follows, see Remark~\ref{remark:equiv}. ISS of \eqref{ws} follows
then by Theorem~\ref{thm:ISS-equiv}. \hfill $\Box$ \vspace{-0.4cm}
\begin{remark}
In comparison to Theorem~1 in \cite{LiN06}, we require
additionally in Theorem~\ref{the:main_traj} that the set $\{f(x,u){:}u{\in}
U{\cap}\epsilon\bar{\B}\}$ is convex for each $x{\in}\chi,
\epsilon{>}0$.
This is due to the fact that we use in our proof that ISS is
equivalent to $0$-input pre-stability and the AG property. This
equivalence requires that the set $\{f(x,u){:}u{\in}
U{\cap}\epsilon\bar{\B}\}$ is convex for each $x{\in}\chi,
\epsilon{>}0$, see \cite{CaT09}. However, we do not exclude that
it might be possible to prove the theorem without this equivalence and to avoid this restriction.
\end{remark}
\vspace{-0.2cm}
\subsection{\hspace{-0.5mm}Small-gain\hspace{-0.2mm} theorems\hspace{-0.2mm} in terms of Lyapunov functions}
\vspace{-0.2cm}
In this section we show how an ISS Lyapunov function for an
interconnection \eqref{ws} can be constructed using the small gain
condition. This allows to apply Proposition~\ref{prop:Lyapunov for
ISS} to deduce ISS of \eqref{ws}.
\vspace{-0.2cm}
\begin{theorem}\label{theorem_iss_lyapunov}
Consider system (\ref{ws}) as interconnection of subsystems
\eqref{is} with $D_i=D, i=1,\dots,n$ and assume that each
subsystem {\it i} has an ISS Lyapunov function $V_i$ satisfying
\eqref{is cond 1 sum}-\eqref{is cond 2.2} with corresponding
ISS-Lyapunov gains. Let the corresponding gain operator
$\Gamma_{\max}$, defined by \eqref{operator_gamma} in terms of
these gains, satisfy \eqref{eq:SGC_max}, then the hybrid system
(\ref{ws}) has an ISS-Lyapunov function given by
\vspace{-0.15cm}
\begin{equation}\label{Lyapunov_function}
V(x)=\max\limits_{i}\sigma^{-1}_i(V_i(x_i))
\end{equation}
\vspace{-0.45cm}\\
where $\sigma_i(r):=\max \{a_i
r,(\Gamma_{\max}(ar))_i,\ldots,(\Gamma_{\max}^{n-1}(ar))_i\}$,
$r{\in}\R_+$ with an arbitrary positive vector
$a{=}(a_1,\dots,a_n)^T.$
In particular, function $V$ satisfies:\\
1) There exist functions $\psi_{1}, \psi_{2} \in \K_{\infty}$ s.t.:
\vspace{-0.15cm}
\begin{equation}\label{ws cond 1}
\psi_{1}(|x|)\leq V(x)\leq  \psi_{2}(|x|)  \mbox{ for any } x\in \chi.
\vspace{-0.15cm}
\end{equation}
2) There exist function $\gamma \in \K$, and continuous, positive
definite functions $\alpha$, $\lambda$ with $\lambda(s)<s$ for all
$s>0$ s.t.:
\vspace{-0.15cm}
\begin{equation}\label{ws cond 2.1}
V(\hspace{-0.05cm}x\hspace{-0.05cm}){\ge}\gamma(\hspace{-0.02cm}|u|\hspace{-0.02cm}){\Rightarrow}\hspace{-0.05cm}{\forall}\zeta{\in}\partial V(x){:}
\hspace{-0.05cm}\left\langle\zeta{,}f\hspace{-0.02cm}(\hspace{-0.02cm}x{,}u\hspace{-0.02cm})\right\rangle\hspace{-0.05cm}{\leq}\hspace{-0.05cm}{-}\alpha(\hspace{-0.02cm}V(\hspace{-0.02cm}x\hspace{-0.02cm})\hspace{-0.02cm}){,}\hspace{-0.02cm}(\hspace{-0.02cm}x{,}u\hspace{-0.02cm})\hspace{-0.02cm}{\in}C{,}
\end{equation}
\vspace{-0.15cm}
\begin{equation}\label{ws cond 2.2}
V(g(x,u)){\leq}\max\{\lambda(V(x),\gamma(|u|)\}, (x,u){\in} D.
\vspace{-0.05cm}
\end{equation}
\end{theorem}
\vspace{-0.35cm}
 \textit{\textbf{Proof.}} First, we establish some
regularity and monotonicity properties of  $\sigma_i$.
Then we apply these properties to show that $V$ constructed
in \eqref{Lyapunov_function} satisfies  \eqref{ws cond 1}-\eqref{ws cond 2.2}.\\
Without loss of generality, the gains $\gamma_{ij}$ can be assumed
to be smooth on $(0,\infty)$, see \cite[Lemma~B.2.1]{Gru02b}. To establish
the properties of $\sigma_i$ consider the map
$Q{:}\R_{+}^n{\to}\R_{+}^n$ defined by
$Q(x){:=}(Q_1(x),{\ldots},Q_n(x))^T$ with $Q_i(s){:=}\max \{s_i,(\Gamma_{\max}(s))_i,{\ldots},(\Gamma_{\max}^{n-1}(s))_i\}$. Note that according to the notation given in
Section~\ref{sec:Preliminaries}, $(\Gamma_{\max}^{p}(s))_i$ is
the $i$th component of the vector $\Gamma_{\max}^{p}(s)$ where
$\Gamma_{\max}^{p}$ is the  $p$-fold composition of
$\Gamma_{\max}$.
By
\cite[Proposition~2.7]{KaJ09}
 the inequality $\Gamma_{\max}(Q(x)){\leq} Q(x)$
holds for all $x{\geq} 0$.
Similarly $\forall\,x>0$ it holds $\Gamma_{\max}(Q(x)){<} Q(x)$.
Fix any positive vector $a>0$ and
consider  $\sigma(r){:=}Q(ar){\in}\R_{+}^n$.
\vspace{-0.2cm}
\begin{equation}\label{Omega_path}
\mbox{Obviously, }\qquad    \Gamma_{\max}(\sigma(r))<\sigma(r),\;\forall r>0
    \vspace{-0.2cm}
\end{equation}
and by the definition of $Q$ it follows that $\sigma_i\in\K_{\infty}$ for all $i=1,\ldots,n$. Furthermore,
$\sigma_i$ satisfy:\\
(i) $\sigma^{-1}_i{\in}\Lip$ on $(0,\infty)$ (as $\gamma_{ij}$ is smooth on $(0,\infty)$);\\
(ii) for every compact set $K\subset(0,\infty)$ there are finite
constants $0<K_1<K_2$ such that for all points of differentiability of
$\sigma^{-1}_i$ we have
\vspace{-0.25cm}
\begin{equation}\label{sigma cond 1}
0<K_1\leq(\sigma^{-1}_i)'(r)\leq K_2, \quad \forall r\in K.
\vspace{-0.25cm}
\end{equation}
In particular, (i) implies that $V$ defined in
\eqref{Lyapunov_function} is locally Lipschitz continuous on
$(0,\infty)$ and (ii) implies the bounded growth of $\sigma_i$ and $\sigma^{-1}_i$ outside the origin. \\
Let us show that such function $V$ satisfies \eqref{ws cond 1}-\eqref{ws cond 2.2}.\\
To this end we define\\
$\psi_{1}(|x|):=\min_{i=1,\ldots,n}\sigma_i^{-1}(\psi_{i1}(c_1|x|))$
and
$\psi_{2}(|x|):=\max_{i=1,\ldots,n}\sigma_i^{-1}(\psi_{i2}(c_2|x|))$
for some suitable positive constants $c_1, c_2$ that depend on the norm
$|\,\cdot\,|$. For example, if $|\,\cdot\,|$ denotes the
infinity norm, then one can take $c_1=c_2=1$. By this choice
the condition (\ref{ws cond 1}) is satisfied.
Define the gain of the whole system by
\vspace{-0.25cm}
\begin{equation}\label{eq:gamma}
\gamma(|u|):=\max\limits_j\{\phi^{-1}(\gamma_j(|u|))\}
\vspace{-0.35cm}
\end{equation}
with $\phi{\in}\K_{\infty}$ such that $\phi(t){\leq}\max\{\max\limits_{j,j\neq i}\gamma_{ij}(\sigma_j(t)),\sigma_i(t)\}$ for all $t\geq 0$. Using \eqref{Omega_path} we obtain for each $i$
\vspace{-0.25cm}
\begin{equation}\label{sigma_max}
\max\{\max\limits_{j,j\neq i}\gamma_{ij}(\sigma_j(r)),\phi(r)\}\leq\sigma_i(r),\forall r>0.
\vspace{-0.25cm}
\end{equation}
Consider any $x\neq 0$, as the case $x=0$ is obvious. Define
\vspace{-0.25cm}
\begin{equation}\label{Mi}
I{:=}\{\hspace{-0.05cm}i{\in}\{1,\ldots,n\}{:}
\sigma^{-1}_i\hspace{-0.05cm}(V_i(x_i)){\geq}\hspace{-0.05cm}\max\limits_{j{,}j\neq i}\sigma^{-1}_j\hspace{-0.05cm}(V_j(x_j))\}
\vspace{-0.35cm}
\end{equation}
i.e. the set of indices $i$ for which the maximum in
\eqref{Lyapunov_function} is attained.\\
 Fix any $i\in I$. If $V(x)\geq\gamma(|u|)$,
then by \eqref{eq:gamma} $\phi(V(x))\geq\gamma_i(|u|)$ and from
\eqref{sigma_max}, \eqref{Mi} we have
\vspace{-0.35cm}
\begin{equation}\notag
\begin{array}{lll}
V_i(x_i)=\sigma_i(V(x))&\geq&\max\{\max\limits_{j,j\neq i}\gamma_{ij}(\sigma_j(V(x))),\phi(V(x))\}\\
&\geq&\max\{\max\limits_{j,j\neq i}\gamma_{ij}(V_j(x_j)),\gamma_i(|u|)\}.
\end{array}
\vspace{-0.3cm}
\end{equation}
To show \eqref{ws cond 2.1} assume $(x,u)\in C$.
As $V$ is obtained through the maximization \eqref{Lyapunov_function},
by 
\cite[p.83]{CLS98} we have that\\
\vspace{-0.15cm}
\begin{equation}\label{thm:conv}
\partial V(x)\subset \mbox{conv}\left\{\bigcup_{i\in I}\partial[\sigma_i^{-1}\circ V_i\circ P_i](x)\right\},
\end{equation}
\vspace{-0.15cm}\\
where $P_i(x):=x_i$. Thus we can use the properties of $\sigma_i$ and $V_i$ to find a bound for $\left\langle \zeta,f(x,u)\right\rangle$, $\zeta\in\partial V$. In particular, by the chain rule  for Lipschitz continuous functions in
\cite[Theorem~2.5]{CLS98}, we have
\vspace{-0.35cm}
\begin{equation}\label{thm:C_i_2}
\partial(\sigma_i^{-1}\hspace{-0.05cm}{\circ}V_i)(x_i){\subset}\{c\zeta_i{:}c{\in}{\partial}\sigma_i^{-1}\hspace{-0.05cm}(V_i(x_i)){,}\zeta_i{\in}\partial V_i(x_i)\},
\vspace{-0.15cm}
\end{equation}
where $c$ is bounded away from zero due to \eqref{sigma cond 1}.
Applying \eqref{is cond 2.1} we obtain for all $\zeta_i\in\partial V_i(x_i)$ that
\vspace{-0.25cm}
\begin{equation}\label{thm:C_i}
\left\langle \zeta_i,f_i(x,u)\right\rangle\leq -\alpha_i(V_i(x_i)).
\vspace{-0.25cm}
\end{equation}
To get a bound independent on $i$ on the right-hand side of \eqref{thm:C_i} define for $\rho{>}0$,
$\widetilde{\alpha}_i(\rho){:=}c_{\rho,i}\alpha_i(\rho){>}0$, where
the constant $c_{\rho,i}{:=}K_1$ with $K_1$ corresponding to the set
$K{:=}\{x_i{\in}\chi_i{:}\rho/2{\leq} |x_i|{\leq} 2\rho\}$ given by
\eqref{sigma cond 1}. And for $r{>}0$ define
$\hat{\alpha}(r){:=}\min\{\widetilde{\alpha}_i(V_i(x_i))|\,|x|{=}r,
V(x)=\sigma_i^{-1}(V_i(x_i))\}{>}0$. Thus using \eqref{thm:C_i_2}-\eqref{thm:C_i}
  we obtain
\vspace{-0.25cm}
\begin{equation}\label{thm:C_i_3}
\left\langle \zeta_i,f_i(x,u)\right\rangle\leq -\hat{\alpha}(|x|) \quad \forall \zeta_i{\in}\partial [\sigma_i^{-1}{\circ}V_i](x_i).
\vspace{-0.25cm}
\end{equation}
The same argument applies for all $i{\in }I$. 
Let us now return to $\zeta{\in}\partial V(x)$. From
\eqref{thm:conv} for any $\zeta{\in}\partial V(x)$ we have that
$\zeta{=}\sum\limits_{i{\in} I}\mu_i c_i \zeta_i$ for suitable
$\mu_i{\geq} 0, \sum_{i{\in} I}\mu_i{=}1$, and with
$\zeta_i{\in}\partial(V_i{\circ} P_i)(x)$ and
$c_i{\in}\partial\sigma_i^{-1}(V_i(x_i))$. Using \eqref{thm:C_i_3}, that
$\left\langle  \zeta_i,f(x,u)\right\rangle\hspace{-0.05cm}{=}
\left\langle  P_i(\zeta_i)\hspace{-0.05cm},\hspace{-0.05cm}f_i(x,u)\right\rangle$ for $\zeta_i{\in}\partial(V_i{\circ} P_i)(x)$ due to the properties of the projection function $P_i$
and that $c_i{>}0$ due to \eqref{sigma cond 1}, it follows that
\vspace{-0.35cm}
\[\left\langle \zeta\hspace{-0.05cm},\hspace{-0.05cm}f(x,u)\right\rangle\hspace{-0.05cm}{=}\hspace{-0.05cm}\sum\limits_{i\in I}\hspace{-0.05cm}\mu_i\hspace{-0.05cm}\left\langle  c_i \zeta_i,f(x,u)\right\rangle\hspace{-0.05cm}{=}\hspace{-0.05cm}
\sum\limits_{i\in I}\hspace{-0.05cm}\mu_i\hspace{-0.05cm}\left\langle  c_i P_i(\zeta_i)\hspace{-0.05cm},\hspace{-0.05cm}f_i(x,u)\right\rangle
\vspace{-0.25cm}\]
\vspace{-0.25cm}
 \[\leq -\sum\limits_{i\in I}\mu_i \hat{\alpha}(|x|)\leq -\hat{\alpha}(|x|)\leq -\hat{\alpha}\circ\psi_2^{-1}\circ V(x).
 \vspace{-0.15cm}\]
 Thus condition \eqref{ws cond 2.1} is satisfied with $\alpha:=\hat{\alpha}\circ\psi_2^{-1}$.\\
To show \eqref{ws cond 2.2} assume now that $(x,u)\in D$. Define
\vspace{-0.25cm}
\begin{equation}\label{lambda_sum}
\lambda(t):=\max\limits_{i,j,i\neq
j}\{\sigma^{-1}_i\circ\gamma_{ij}\circ\sigma_j(t),\sigma^{-1}_i\circ\lambda_i\circ\sigma_i(t)\},\;
t>0.
\vspace{-0.05cm}
\end{equation}
Note that $\sigma^{-1}_i\circ\gamma_{ij}\circ\sigma_j(t)<\sigma^{-1}_i\circ\sigma_i(t)=t$
from (\ref{Omega_path}) and $\sigma^{-1}_i\circ\lambda_i\circ\sigma_i(t)<\sigma^{-1}_i\circ\sigma_i(t)=t$ for all $t>0$ as $\lambda_i(t){<}t$. Thus $\lambda(t){<}t,$ $\forall\,t{>}0$. Let us show that such $\lambda$ satisfies \eqref{ws cond 2.2}.
Condition (\ref{is cond 2.2}) for ISS-Lyapunov function
of subsystem $i$,
the jump behaviour \eqref{new-jumps} and the assumption $D_i=D$, $i=\{1,\ldots,n\}$
 imply for $(x,u)\in D$
\vspace{-0.2cm}\\
\begin{equation*}
\begin{split}
V(g(x,u)){=}\max\limits_i\sigma_i^{-1}{\circ} V_i(g_i(x,u))\\
{\leq}\max\limits_{i{,}j{,}i{\neq} j}\{\sigma^{-1}_i{\circ}\lambda_i(V_i(x_i)),
\sigma^{-1}_i{\circ}\gamma_{ij}(V_j(x_j)),\sigma^{-1}_i{\circ}\gamma_i(|u|)\}\\
=\max\limits_{i,j,i\neq j}\{\sigma^{-1}_i\circ \lambda_i\circ\sigma_i\circ\sigma^{-1}_i(V_i(x_i)),\\
\sigma^{-1}_i{\circ}\gamma_{ij}{\circ}\sigma_j{\circ}\sigma^{-1}_j(V_j(x_j)),\sigma^{-1}_i{\circ}\gamma_i(|u|)\}\\
\leq\max\{\lambda(V(x)),\gamma(|u|).
\end{split}
\end{equation*}
Thus \eqref{ws cond 2.2} is also satisfied and hence $V$ is an ISS-Lyapunov function of the
network \eqref{ws}. \hfill $\Box$\\
\vspace{-0.6cm}
\vspace{-1mm}
\section{Conclusions}\label{sec:concl}
\vspace{-3mm} We have shown that a large scale interconnection of
ISS hybrid systems is again ISS if the small gain condition
 is satisfied. The results are provided in terms
of trajectories and Lyapunov functions. Moreover an explicit
construction of an ISS-Lyapunov function is given.
 These results
extend the corresponding known theorems from \cite{NeT08} to the
case of interconnection of more than two hybrid systems and \cite{KaJ09} for
general type of hybrid systems. However, our results are restricted to interconnections with a common jump set of subsystems.
\vspace{-0.35cm}

\begin{ack}
\vspace{-0.2cm} The authors are thankful to the anonymous reviewers
for their useful remarks that led to an improved presentation of the
paper. M. Kosmykov was supported by the Volkswagen Foundation
(Project I/82684). S. Dashkovskiy was supported by the DFG as a part
of the CRC~637.
\end{ack}
\appendix
\vspace{-0.3cm}
\section{Equivalent definition of an ISS-Lyapunov function}\label{section:equivalent_formulations}
\vspace{-0.3cm}
Here we show equivalence between the definition of an ISS-Lyapunov function
used in \cite{CaT09} and Definition~\ref{def:lyap_v_i}.\\
Consider a function $W:\chi\to \R_+$, $W\in\Lip$
that satisfies the following properties for (\ref{ws})\\
1) There exist functions $\bar{\psi}_{1},\bar{\psi}_{2} \in \K_{\infty}$ such that:
\vspace{-0.2cm}
\begin{equation}\label{ws cond 1'}
\bar{\psi}_{1}(|x|)\leq W(x)\leq\bar{\psi}_{2}(|x|)  \mbox{ for any } x\in \chi.
\vspace{-0.2cm}
\end{equation}
2) There exist function $\bar{\gamma} \in \K$, continuous, positive
definite function $\bar{\alpha}_1$ and function
$\bar{\alpha}_2\in\K_{\infty}$ such that:
\vspace{-0.2cm}
\begin{equation}\label{ws cond 2.1'}
|x|{\ge}\bar{\gamma}(|u|){\Rightarrow}{\forall}\zeta{\in}\partial
W(x){:}\left\langle\zeta{,}f(x{,}u)\right\rangle{\le}
{-}\bar{\alpha}_1(|x|){,} (x{,}u){\in}C{,}
\end{equation}
\begin{equation}\label{ws cond 2.2'}
|x|\geq \bar{\gamma}(|u|)\Rightarrow
W(g(x,u))-W(x){\leq}-\bar{\alpha}_2(|x|){,} (x,u){\in} D.
\end{equation}
In \cite{CaT09} the conditions \eqref{ws cond 1'}-\eqref{ws cond
2.2'} with $\bar\alpha_1\in\K_\infty$ were used to define an
ISS-Lyapunov function for \eqref{ws} and it was shown that
existence of such (smooth) function $W$ implies that \eqref{ws} is
ISS. This proof does not change if $\bar\alpha_1$ is continuous and
positive definite only.
\vspace{-0.3cm}
\begin{proposition}\label{prop:lyap_equiv}
System \eqref{ws} has an ISS-Lyapunov function $V$ satisfying
\eqref{ws cond 1}-\eqref{ws cond 2.2} if and only if there exists
$W\in\Lip$ satisfying \eqref{ws cond 1'}-\eqref{ws cond 2.2'}.
\end{proposition}
\vspace{-0.25cm}
\textit{\textbf{Proof.}}
"$\Rightarrow$" \quad
Let $V$ satisfy (\ref{ws cond 1})-(\ref{ws cond 2.2}).
We can always majorize a continuous, positive definite function
$\lambda<\id$ from \eqref{ws cond 2.2} by a
function $\rho\in\K_{\infty}$ such that
$\lambda(r)\leq\rho(r)<r$ for $r>0$, for example,
$\rho:=\frac{1}{2}(\max_{[0,r]}\lambda+\id).$ Then for
$(x,u)\in D$ from \eqref{ws cond 2.2}  we have
\vspace{-0.3cm}
\begin{equation}\label{proof:ws cond 2.2'''}
V\hspace{-0.025cm}(g(x{,}u))\hspace{-0.025cm}{\leq}\hspace{-0.025cm}\max\{\lambda(V(x))\hspace{-0.025cm}{,}\hspace{-0.025cm}\gamma(|u|)\}{\leq}\max\{\rho(\hspace{-0.025cm}V(x)\hspace{-0.025cm}){,}\gamma(\hspace{-0.025cm}|u|\hspace{-0.025cm})\}.
\end{equation}
Define
\vspace{-0.05cm}
\begin{equation}\label{proof:lyap_equiv_bar_gamma}
\bar{\gamma}(|u|):=\psi_1^{-1}\circ\rho^{-1}\circ\gamma.
\end{equation}
\vspace{-0.2cm}
If $|x|{\geq}\bar{\gamma}(|u|)$, then
$\rho(V(x))\geq\gamma(|u|)$ and using
\eqref{proof:ws cond 2.2'''}
\vspace{-0.2cm}
\begin{equation}\notag
V(g(x,u))\hspace{-0.025cm}{\leq}\hspace{-0.025cm}\max\{\rho(V\hspace{-0.025cm}(x)){,}\gamma(|u|)\}{=}\rho(V\hspace{-0.025cm}(x)){=}V\hspace{-0.025cm}(x){-}\widetilde{\alpha}(x)
\end{equation}
\vspace{-0.5cm}
\begin{equation}\label{proof:ws cond 2.2''}
\Rightarrow V(g(x,u))-V(x)\leq-\hat{\alpha}(|x|),
\end{equation}
with $\hat{\alpha}(r){:=}\min\limits_{|s|{=}r}\widetilde{\alpha}(s)$ that is a continuous,
positive definite function, where $\widetilde{\alpha}(s){:=}V(s){-}\rho(V(s)){\geq} 0$. From \eqref{proof:ws cond 2.2''} and \cite[Lemma~2.8]{JiW01b},  ${\exists}\bar\rho,\bar{\alpha}_2{\in}\K_{\infty}$ such that $W{:=}\bar\rho{\circ} V$ satisfies \eqref{ws
cond 2.2'} with $\bar{\gamma}(|u|)$ defined in
\eqref{proof:lyap_equiv_bar_gamma}.
As $V(x)$ satisfies \eqref{ws cond 2.1}, then $W$ satisfies \eqref{ws cond 2.1'} with $\bar{\alpha}_1{:=}\widetilde{\rho}{\cdot}\alpha$ that is continuous, positive definite function, where $\widetilde{\rho}{\in}\partial\bar\rho(V(x))$.\\
Thus function  $W$ satisfies \eqref{ws cond 1'}-\eqref{ws cond 2.2'} with $\bar{\psi}_i{:=}\bar\rho{\circ} \psi_i$. 
"$\Leftarrow$" \quad
Assume now that function $W$ satisfies \eqref{ws cond 1'}-\eqref{ws cond 2.2'} 
and define $V:=W$, $\psi_1:=\bar{\psi}_1$ and $\psi_2:=\bar{\psi}_2$. Then condition \eqref{ws cond 1} is satisfied.
Let
\vspace{-0.25cm}
\begin{equation}\label{proof:lyap_equiv_gamma}
\gamma(|u|):=\bar{\psi}_2\circ\bar{\gamma}(|u|).
\vspace{-0.25cm}
\end{equation}
Consider $V(x)\geq\gamma(|u|)$. Then from \eqref{proof:lyap_equiv_gamma}, \eqref{ws cond 1'} 
$|x|\geq\bar{\gamma}(|u|)$. From \eqref{ws cond 1'}-\eqref{ws cond 2.1'} for all $(x,u)\in C$
\vspace{-0.25cm}
\begin{equation*}
\forall\zeta\in\partial V(x){:}\left\langle \zeta,f(x,u)\right\rangle{\leq}{-}\bar{\alpha}_1(|x|)
\leq{-}\bar{\alpha}_1\circ\bar{\psi}_2^{-1}(V(x)).
\vspace{-0.3cm}
\end{equation*}
Thus $V$ satisfies \eqref{ws cond 2.1} with $\alpha:=\bar{\alpha}_1\circ\bar{\psi}_2^{-1}$.\\
By \cite[Lemma~B.1]{JiW01} for any
$\bar{\alpha}_2{\circ}\bar{\psi}_2^{-1}{\in}\K_{\infty}$
$\exists\widetilde{\alpha}{\in}\K_{\infty}$ such that
$\widetilde{\alpha}{\leq}\bar{\alpha}_2{\circ}\bar{\psi}_2^{-1}$ and
$\id{-}\widetilde{\alpha}{\in}\K$. For $V(x)>\gamma(|u|)$ from \eqref{ws cond 1'} and
\eqref{ws cond 2.2'}
\vspace{-0.35cm}
\begin{equation*}
\begin{array}{l}
V(g(x,u)){\leq} V(x){-}\bar{\alpha}_2(|x|){\leq} V(x){-}\bar{\alpha}_2\circ\bar{\psi}_2^{-1}(V(x))\\
{\leq}V(x){-}\widetilde{\alpha}(V(x)){=}(\id-\widetilde{\alpha})(V(x)){=}\lambda(V(x)),
\end{array}
\vspace{-0.35cm}
\end{equation*}
where $\lambda:=\id-\widetilde{\alpha}$.\\
Consider now
$(x,u){\in} D$ such that $V(x){\leq}\gamma(|u|)$ and define
$\A(|u|){:=}\{(x,u){\in} D:V(x)\leq\gamma(|u|)\}$.\\
 Let us take now $\hat{\gamma}(|u|){:=}\max\limits_{(x,u){\in} \A(|u|)}V(g(x,u))$. Then $V(g(x,u)){\leq} \hat{\gamma}(|u|)$.
  Note that $\hat{\gamma}(0){=}0$ as $V(x){\geq} 0{=}\gamma(0)$. Furthermore, as function $V$ is nonnegative and $V{\in}\Lip$ and function $g$ is continuous,
  function $\hat{\gamma}{\in}\Lip$ is nonnegative.
 We can always majorize such function $\hat{\gamma}$  by a function $\check{\gamma}{\in}\K$ such that $\hat{\gamma}{\leq}\check{\gamma}$. Thus for $(x,u){\in} D$ we obtain that $V(g(x,u)){\leq}\max\{\check{\gamma}(|u|),\lambda(V(x))\}$ and condition \eqref{ws cond 2.2} is satisfied with $V{:=}W$ and $\widetilde{\gamma}{:=}\max\{\check{\gamma},\gamma\}$.
\hfill $\Box$\\
\vspace{-0.4cm}
\vspace{-3mm}
\section{Proofs of Theorem~\ref{sgc_gs} and Theorem~\ref{sgc_ag}}
\vspace{-0.25cm}
We need first the following auxiliary lemmas.
\vspace{-0.25cm}
\begin{lemma}\cite[Theorem~2.5.4.]{Rue07}\label{lemma_w_bound_max}
Let the operator $\Gamma_{\max}$ defined in
\eqref{operator_gamma} satisfy \eqref{eq:SGC_max}. Then there exists $\phi \in
\K_{\infty}$ s.t. for all $w, v \in {\mathbb{R}}_{+}^n$,
\vspace{-0.35cm}
\begin{equation*}
w\leq\max\{\Gamma_{\max}(w),v\}:=\left(\begin{array}{c}
\max\{\max\limits_{j,j\neq 1}\gamma_{1j}(w_j),v_1\}\\
\vdots\\
\max\{\max\limits_{j,j\neq n}\gamma_{nj}(w_j),v_n\}
\end{array}\right)
\vspace{-0.25cm}
\end{equation*}
implies $\|w\|\leq\phi(\|v\|)$.
\end{lemma}
\vspace{-0.25cm}
\begin{lemma}\label{lemma_limsup_ag}
Let $s:\dom \,s\to\R_{+}^n$
be continuous between the jumps and bounded with unbounded $\dom \,s$. Then\\
$\limsup\limits_{(t,k)\in\doms\, s,t+k\to\infty}\hspace{-0.25cm}s(t,k)=\limsup\limits_{t+k\to\infty}\|s_{[(t/2,k/2),\lim\limits_{\tau+j\to\infty}(\tau,j)]}\|$.
\vspace{-0.25cm}
\end{lemma}
\vspace{-0.25cm}
\textit{\textbf{Proof.}}
The proof goes along the lines of the proof of a similar result for continuous systems in Lemma~3.2 in \cite{DRW07} but instead of time $t$ we consider the points $(t,k)$ of the time domain.
\hfill $\Box$
\vspace{-0.25cm}
\subsection{Proof of Theorem~\ref{sgc_gs}}\label{appendix:sgc_gs}
\vspace{-3mm} Let us take the supremum over $(\tau,l)\leq (t,k)$
on both sides of \eqref{iGS_max}
\vspace{-0.4cm}
\begin{equation}\label{sup_gs}
\begin{split}
{\|x_{i(t,k)}\|}_{(\bar{\tau}{,}\bar{l})}{\leq}
\max\{\sigma_{i}(|x_{i0}|){,}\\
\max\limits_{j{,}j{\neq} i}\hat{\gamma}_{ij}({\|x_{j(t,k)}\|}_{(\bar{\tau}{,}\bar{l})}),\hat{\gamma}_{i}({\|u\|}_{(\bar{\tau}{,}\bar{l})})\},
\end{split}
\vspace{-0.45cm}
\end{equation}
\vspace{-0.45cm}\\
where $(\bar{\tau},\bar{l}):=\max\limits_{(\tau,l)\in \doms\,x} (\tau, l)$, i.e. the maximum element of $\dom\,x$.
\\
Let $\Gamma{:=}(\hat{\gamma}_{ij})_{n\times n}$,
$w{:=}\left({\|x_{1(t,k)}\|}_{(\bar{\tau},\bar{l})},\dots,{\|x_{n(t,k)}\|}_{(\bar{\tau},\bar{l})}\right)^T$,
\vspace{-0.15cm}
\begin{equation*}
\begin{array}{lll}
\quad\quad\quad v&:=&\left(\begin{array}{c}
\max\{\sigma_{1}(|x_{10}|),\hat{\gamma}_{1}({\|u\|}_{(\bar{\tau},\bar{l})})\}\\
\vdots\\
\max\{\sigma_{n}(|x_{n0}|),\hat{\gamma}_{n}({\|u\|}_{(\bar{\tau},\bar{l})})\}
\end{array}\right)\\
\quad\quad\quad &=&\max\{\sigma(|x_0|),\hat{\gamma}({\|u\|}_{(\bar{\tau},\bar{l})})\},
\end{array}
\vspace{-0.15cm}
\end{equation*}
where $\sigma,\hat{\gamma}\in\K_{\infty}$. From \eqref{sup_gs} we
obtain $w\leq\max\{\Gamma_{\max}(w),v\}$. Then by
Lemma~\ref{lemma_w_bound_max} $\exists\,\phi\in\K_{\infty}$ s.t.
\vspace{-0.35cm}
\begin{equation}\label{gs bound}
\begin{array}{lll}
|x(t,k)|&{\leq}&{\|x_{(t,k)}\|}_{(\bar{\tau},\bar{l})}{\leq}\phi(\|\max\{\sigma(|x_0|),\hat{\gamma}({\|u\|}_{(\bar{\tau},\bar{l})})\}\|)\\
&{\leq}&\max\{\phi(\|\sigma(|x_0|)\|),\phi(\|\hat{\gamma}({\|u\|}_{(\bar{\tau},\bar{l})})\|)\}.
\end{array}
\end{equation}
for all $(t,k)\in\dom\,x$. Hence for
every initial condition and essentially bounded input $u$ the
solution of the system (\ref{ws})  exists 
and is bounded, since the right-hand side of (\ref{gs
bound}) does not depend on $t,k$. From the last line in (\ref{gs bound}) the estimate \eqref{ws_GS_max} for the pre-GS follows. \hfill $\Box$
\vspace{-0.25cm}
\subsection{Proof of Theorem~\ref{sgc_ag}}\label{appendix:sgc_ag}
\vspace{-3mm} 
Let $(\tau,l)$ be an arbitrary initial point of the time domain. From the
definition of the AG property we have
\vspace{-0.35cm}
\begin{equation*}
\begin{split}
\limsup\limits_{(t,k)\in\doms \, x_i,t+k\to\infty}|x_i(t,k)| \\
\leq\max\{\max\limits_{j,j\neq
i}\widetilde{\gamma}_{ij}({\|x_{j[(\tau,l),\lim\limits_{\bar{\tau}+\bar{l}\to\infty}(\bar{\tau},\bar{l})]}}\|_{\infty}),
\widetilde{\gamma}_i({\|u\|}_{\infty})\}.
\end{split}
\vspace{-0.45cm}
\end{equation*}
\vspace{-0.45cm}\\
Then from Lemma~3.6 in \cite{CaT09} it follows that
\vspace{-0.35cm}
\begin{equation}\label{ag_proof_sum_2}
\begin{split}
\limsup\limits_{(t,k)\in\doms \, x_i,t+k\to\infty}\hspace{-0.4cm}|x_i(t,k)| \\
{\leq}\hspace{-0.05cm}\max\{\max\limits_{j,j{\neq}
i}\hspace{-0.05cm}\widetilde{\gamma}_{ij}(\limsup\limits_{\tau{+}l {\to}
 \infty}({\|x_{j[(\tau,l){,}\hspace{-0.2cm}\lim\limits_{\bar{\tau}{+}\bar{l}\to\infty}\hspace{-0.2cm}(\bar{\tau}{,}\bar{l})]}}\|_{\infty})){,}
\hspace{-0.05cm}\widetilde{\gamma}_i({\|u\|}_{\infty})\}.
\end{split}\end{equation}
\vspace{-0.4cm}\\
Since all solutions of (\ref{is}) are bounded
the following holds by Lemma~\ref{lemma_limsup_ag}:\\
\vspace{-0.35cm}
\[\limsup\limits_{
\substack{(t,k)\in\doms \, x_i{,}\\
            t{+}k\to\infty}}
{\hspace{-0.3cm}|x_i(t{,}k)|}
{=}\hspace{-0.05cm}\limsup\limits_{\tau{+}l {\to}
\infty}(\|x_{i[(\tau,l){,}\hspace{-0.2cm}\lim\limits_{\bar{\tau}{+}\bar{l}\to\infty}\hspace{-0.1cm}(\bar{\tau},\bar{l})]}\|_{\infty}\hspace{-0.05cm}){=:}l_i(x_i){.}
\vspace{-0.10cm}\]
By this property from (\ref{ag_proof_sum_2}) follows
\vspace{-0.35cm}
\begin{equation}\notag
l_i(x_i)\leq
\max\{\max\limits_{j,j\neq i}\widetilde{\gamma}_{ij}(l_j(x_j)),\widetilde{\gamma}_{i}({\|u\|}_{\infty})\}.
\vspace{-0.35cm}
\end{equation}
Using Lemma \ref{lemma_w_bound_max} for $\Gamma=(\widetilde{\gamma}_{ij})_{n\times n}$, $w_i=l_i(x_i)$ and $v_i:=\widetilde{\gamma}_{i}({\|u\|}_{\infty})$ we conclude
\vspace{-0.35cm}
\begin{equation}\label{bound ag}
\limsup\limits_{(t,k)\in\doms \, x,t+k\to\infty}|x(t,k)|\leq\phi({\|u\|}_{\infty})
\vspace{-0.35cm}
\end{equation}
for some $\phi{\in}\K$, which is the desired AG property. \hfill
$\Box$ \vspace{-5mm}
\bibliographystyle{plain}        
\bibliography{literatur}
\end{document}